\documentclass[12pt,twoside]{amsart}
\usepackage{amssymb,amsmath,amsthm, amscd, enumerate, mathrsfs}
\usepackage{bigdelim}
\usepackage{color} 
\usepackage{multirow}
\usepackage[all]{xy} 
\usepackage{comment}
\usepackage{graphicx, hhline}
\usepackage[dvipdfmx]{hyperref}

\title{On vanishing theorems for analytic spaces}
\author{Osamu Fujino}
\date{2023/10/15, version 0.02}
\address{Department of 
Mathematics, Graduate School of Science, 
Kyoto University, Kyoto 606-8502, Japan}
\email{fujino@math.kyoto-u.ac.jp}

\DeclareMathOperator{\Ass}{Ass}
\DeclareMathOperator{\Pic}{Pic}
\DeclareMathOperator{\mult}{mult}
\newtheorem{thm}{Theorem}[section]

\theoremstyle{definition}
\newtheorem*{ack}{Acknowledgments}  
\newtheorem{say}[thm]{}
\begin{document}

\maketitle

\begin{abstract} 
This is an announcement of our 
new vanishing theorems for projective 
morphisms between complex analytic spaces. 
We established a complex analytic 
generalization of Koll\'ar's torsion-freeness 
and vanishing theorem for analytic simple normal 
crossing pairs. Although our results may look artificial, 
they have already played a crucial role for the study 
of minimal models in the complex analytic setting. 
\end{abstract}

\section{Introduction}\label{b-sec1}
This is a short announcement of our 
new vanishing theorems for projective 
morphisms between complex analytic spaces.
All the details will be published in \cite{fujino-vanishing-ne}. 

In his monumental paper \cite{kollar1}, 
Koll\'ar generalized the Kodaira 
vanishing theorem for complex projective varieties. His results consist 
of injectivity, torsion-free, and 
vanishing theorems. We have already had a 
powerful generalization of Koll\'ar's 
package for {\em{reducible}} algebraic varieties 
(for the details, see, for example, 
\cite[Chapter 5]{fujino-foundations}), 
which plays a crucial role for the study of 
log canonical pairs, semi-log canonical 
pairs, and quasi-log schemes in the theory of 
minimal models of algebraic varieties 
(see \cite{fujino-fundamental}, 
\cite{fujino-fundamental-slc}, 
\cite[Chapter 6]{fujino-foundations}, and 
so on). 
Hence it was highly 
desirable to establish an analytic generalization (see 
\cite[Remark 5.8.3]{fujino-foundations} and \cite[1.10]{fujino-minimal}). 
Roughly speaking, from the Hodge theoretic viewpoint, 
Koll\'ar's original result in \cite{kollar1} 
is {\em{pure}} 
and the generalization in \cite[Chapter 5]{fujino-foundations} is 
{\em{mixed}}. 
Recently, in \cite{fujino-vanishing-ne}, we established 
an appropriate generalization of \cite[Chapter 5]{fujino-foundations} 
for projective morphisms of complex analytic spaces. 
By this new generalization, 
we can translate the results in \cite{fujino-fundamental}, 
\cite{fujino-fundamental-slc}, and \cite[Chapter 6]{fujino-foundations} 
into the ones for projective morphisms between complex 
analytic spaces (see \cite{fujino-cone-contraction} and 
\cite{fujino-quasi-log-analytic}). 
More precisely, in \cite{fujino-cone-contraction}, 
we proved the 
cone and contraction theorem of normal pairs for 
projective morphisms between complex analytic spaces as an application 
of \cite{fujino-vanishing-ne}. 
Then, in \cite{fujino-quasi-log-analytic}, 
we discussed quasi-log structures 
for complex analytic spaces. 
We have already established the theory of minimal models 
for projective morphisms of complex analytic spaces 
with mild singularities in \cite{fujino-minimal}, which 
is an analytic generalization of the 
great work of Birkar--Cascini--Hacon--M\textsuperscript{c}Kernan. 
We note that \cite{fujino-minimal} does not need our new vanishing theorems. 
The Kawamata--Viehweg vanishing theorem for projective morphisms 
of complex analytic spaces is sufficient for \cite{fujino-minimal}. 
Finally, we recommend the reader who is interested in 
vanishing theorems and the minimal model program to 
see \cite[Chapter 3]{fujino-foundations}. 

In this paper, 
every complex analytic space is assumed to be 
{\em{Hausdorff}} and {\em{second-countable}}. 
We will freely use the standard notation in \cite{fujino-fundamental}, 
\cite{fujino-foundations}, \cite{fujino-minimal}, and so on. 
Let us prepare various definitions in order to 
explain our new vanishing theorems. 

\begin{say}[Analytic globally embedded simple normal crossing 
pairs]\label{b-say1.1}
Let $X$ be a simple normal crossing divisor 
on a smooth complex analytic space $M$ and 
let $B$ be an $\mathbb R$-divisor on $M$ such that 
the support of $B+X$ is a simple normal crossing divisor on $M$ and 
that $B$ and $X$ have no common irreducible components. 
Then we put $D:=B|_X$ and 
consider the pair $(X, D)$. 
The pair $(X, D)$ is called an {\em{analytic globally embedded simple 
normal crossing pair}}. 
\end{say}

Analytic globally embedded simple normal crossing pairs naturally appear 
when we use the resolution of singularities. 

\begin{say}[Analytic simple normal crossing pairs]\label{b-say1.2}
If the pair $(X, D)$ is locally isomorphic to an analytic 
globally embedded 
simple normal crossing pair and the irreducible 
components of $X$ and $D$ are all smooth, 
then $(X, D)$ is called an {\em{analytic simple normal crossing 
pair}}. 
When $(X, D)$ is an analytic simple normal crossing pair, $X$ has an invertible 
dualizing sheaf $\omega_X$.
\end{say}
\begin{say}[Strata]\label{b-say1.3}
Let $(X, D)$ be an analytic simple normal crossing pair. 
Let $\nu\colon X^\nu\to X$ be the normalization. 
Let $\Theta$ be the union of $\nu^{-1}_*D$ and the 
inverse image of the singular locus of $X$. 
If $W$ is an irreducible component of $X$ 
or the $\nu$-image 
of some log canonical center of $(X^\nu, \Theta)$, 
then $W$ is called a {\em{stratum}} of $(X, D)$. 
Note that $X^\nu$ is smooth and the support of $\Theta$ 
is a simple normal crossing divisor on $X^\nu$. 
We also note that if the coefficients of $D$ are in $[0, 1]$ then 
$(X^\nu, \Theta)$ is log canonical. 
\end{say}

In the theory of minimal models, the notion of $\mathbb R$-line bundles 
is indispensable. 

\begin{say}[$\mathbb R$-line bundles and $\mathbb Q$-line bundles]
\label{b-say1.4}
Let $X$ be a complex analytic space and let $\Pic(X)$ 
be the group of line bundles on $X$, 
that is, the {\em{Picard group}} of $X$. An element 
of $\Pic(X)\otimes _{\mathbb Z}\mathbb R$ 
(resp.~$\Pic(X)\otimes _{\mathbb Z}\mathbb Q$) is called 
an {\em{$\mathbb R$-line bundle}} 
(resp.~a {\em{$\mathbb Q$-line bundle}}) on $X$. 
As usual, in this paper, we write the group law of 
$\Pic(X)\otimes _{\mathbb Z}\mathbb R$ additively 
for simplicity of notation. 
\end{say}

We need Siu's theorem to state our result. 

\begin{say}[{Associated subvarieties, see \cite{siu}}]\label{b-say1.5}
Let $\mathcal F$ be a coherent sheaf on a complex 
analytic space $X$. 
Then there exists a locally finite family $\{Y_i\}_{i\in I}$ 
of complex analytic subvarieties of $X$ such that 
\begin{equation*}
\Ass _{\mathcal O_{X,x}}(\mathcal F_x)=\{\mathfrak{p}_{x, 1}, 
\ldots, \mathfrak{p}_{x, r(x)}\}
\end{equation*}
holds for every point $x\in X$, where 
$\mathfrak{p}_{x, 1}, 
\ldots, \mathfrak{p}_{x, r(x)}$ are the prime ideals 
of $\mathcal O_{X, x}$ associated to the irreducible components 
of the germs $Y_{i, x}$ of $Y_i$ at $x$ with $x\in Y_i$. 
We note that each $Y_i$ is called an {\em{associated subvariety}} 
of $\mathcal F$. 
\end{say}

The following theorem is the main result of \cite{fujino-vanishing-ne}, 
which is obviously an analytic generalization of 
\cite[Theorem 5.6.2]{fujino-foundations}. 

\begin{thm}
[{Main theorem, \cite[Theorem 1.1]{fujino-vanishing-ne}}]\label{b-thm1.6}
Let $(X, \Delta)$ be an analytic simple 
normal crossing pair such that the coefficients of 
$\Delta$ are in $[0, 1]$. 
Let $f\colon X\to Y$ be a projective morphism 
to a complex analytic space $Y$ and let $\mathcal L$ 
be a line bundle on $X$. 
Let $q$ be an arbitrary nonnegative integer. 
Then we have the following properties. 
\begin{itemize}
\item[(i)] $($Strict support condition$)$.~If 
$\mathcal L-(\omega_X+\Delta)$ is $f$-semiample,  
then every 
associated subvariety of $R^qf_*\mathcal L$ is the $f$-image 
of some stratum of $(X, \Delta)$. 
\item[(ii)] $($Vanishing theorem$)$.~If 
$\mathcal L-(\omega_X+\Delta)\sim _{\mathbb R} f^*\mathcal H
$ holds 
for some $\pi$-ample 
$\mathbb R$-line bundle $\mathcal H$ on $Y$, where 
$\pi\colon Y\to Z$ is a 
projective morphism to a complex analytic space 
$Z$, then we have 
$
R^p\pi_*R^qf_*\mathcal L=0
$
for every $p>0$. 
\end{itemize} 
\end{thm}

Since we treat complex analytic spaces, we have to 
be careful about some basic definitions. 

\begin{say}[Globally $\mathbb R$-Cartier 
divisors]\label{b-say1.7}
In Theorem \ref{b-thm1.6}, we 
always implicitly 
assume that $\Delta$ is {\em{globally $\mathbb R$-Cartier}}, 
that is, $\Delta$ is a finite $\mathbb R$-linear combination 
of Cartier divisors on $X$. We note that 
if the number of the irreducible components of the 
support of 
$\Delta$ is finite then $\Delta$ is globally $\mathbb R$-Cartier. 
This condition is harmless to 
applications 
because the restriction of $\Delta$ to a relatively compact 
open subset of $X$ has only 
finitely many irreducible components in its support. 
Under the assumption that 
$\Delta$ is globally $\mathbb R$-Cartier, 
we can obtain an $\mathbb R$-line bundle 
$\mathcal N$ naturally associated to $\mathcal L-(\omega_X+\Delta)$, 
which is a hybrid of line bundles $\mathcal L$ and 
$\omega_X$ and a globally $\mathbb R$-Cartier divisor 
$\Delta$. 
The assumption in Theorem \ref{b-thm1.6} (i) means that 
$\mathcal N$ is a finite positive $\mathbb R$-linear 
combination of $f$-semiample line bundles on $X$. 
The assumption in Theorem \ref{b-thm1.6} (ii) 
means that $\mathcal N=f^*\mathcal H$ holds 
in $\Pic(X)\otimes _{\mathbb Z} \mathbb R$.  
\end{say}

\begin{say}\label{b-say1.8}
Let $X=\mathbb C$ and let $\{P_n\}_{n=1}^\infty$ be a 
set of mutually distinct discrete points of $X$. 
Then $\Delta=\sum _{n=1}^\infty \frac{1}{n} P_n$ is a 
$\mathbb Q$-Cartier $\mathbb Q$-divisor on $X$. 
However, it is not a finite $\mathbb R$-linear 
combination of Cartier divisors on $X$. 
Hence it is not a globally $\mathbb R$-Cartier divisor. 
\end{say}

Theorem \ref{b-thm1.6} (ii) can be generalized as follows. 
It is an analytic generalization of \cite[Theorem 5.7.3]{fujino-foundations}. 

\begin{thm}[{Vanishing theorem of Reid--Fukuda type, 
\cite[Theorem 1.2]{fujino-vanishing-ne}}]\label{b-thm1.9}
Let $(X, \Delta)$ be an analytic simple 
normal crossing pair such that the coefficients of 
$\Delta$ are in $[0, 1]$. 
Let $f\colon X\to Y$ and $\pi\colon Y\to Z$ be projective morphisms 
between complex analytic spaces and let $\mathcal L$ 
be a line bundle on $X$. 
If $\mathcal L-(\omega_X+\Delta)\sim _{\mathbb R} f^*\mathcal H$ 
holds such that $\mathcal H$ is an $\mathbb R$-line bundle, 
which is nef and 
log big over $Z$ with respect to $f\colon (X, \Delta)\to Y$, then  
$R^p\pi_*R^qf_*\mathcal L=0$ holds for every $p>0$ and every $q$. 
\end{thm}

The definition of nef and log big line bundles 
in Theorem \ref{b-thm1.9} is as follows. 

\begin{say}[Nef and log big line bundles]\label{b-say1.10}
In Theorem \ref{b-thm1.9}, we 
note that $\mathcal H$ is said to be {\em{nef and 
log big over $Z$ with 
respect to $f\colon (X, \Delta)\to Y$}} if 
$\mathcal H\cdot C\geq 0$ holds for 
every projective integral curve $C$ on $Y$ such that 
$\pi(C)$ is a point and 
$\mathcal H|_{f(W)}$ can be written as a finite 
positive $\mathbb R$-linear combination of $\pi$-big line bundles on 
$f(W)$ for every stratum $W$ of $(X, \Delta)$.
\end{say}

It is more or less well known that Theorem \ref{b-thm1.6} 
follows from Theorem \ref{b-thm1.11} (i) and (ii) below and 
that Theorem \ref{b-thm1.11} (iii) is an easy consequence of 
Theorem \ref{b-thm1.11} (i) and (ii). 
Hence all we have to do is to establish Theorem \ref{b-thm1.11} 
(i) and (ii). We will prove them in Section \ref{b-sec2}. 

\begin{thm}[Koll\'ar's package for 
analytic simple normal crossing pairs]\label{b-thm1.11}
Let $(X, D)$ be an analytic simple normal crossing pair 
such that $D$ is reduced and 
let $f\colon X\to Y$ be a projective  
morphism of complex analytic spaces. 
Then we have the following properties. 
\begin{itemize}
\item[(i)] $($Strict support condition$)$.~Every 
associated subvariety of $R^qf_*\omega_X(D)$ 
is the $f$-image of some stratum of $(X, D)$ for every $q$. 
\item[(ii)] $($Vanishing theorem$)$.~Let $\pi\colon Y\to Z$ be 
a projective morphism between complex analytic 
spaces and let $\mathcal A$ be a $\pi$-ample 
line bundle on $Y$. 
Then 
\begin{equation*}
R^p\pi_*\left(\mathcal A\otimes R^qf_*\omega_X(D)\right)=0
\end{equation*} 
holds  
for every $p>0$ and every $q$. 
\item[(iii)] $($Injectivity theorem$)$.~Let 
$\mathcal L$ be an $f$-semiample 
line bundle on $X$. Let $s$ be a nonzero element of $H^0(X, 
\mathcal L^{\otimes k})$ for some nonnegative integer $k$ such that 
the zero locus of $s$ does not contain any strata of $(X, D)$. 
Then, for every $q$, 
the map 
\begin{equation*}
\begin{split}
\times s&\colon R^qf_*\left(\omega_X(D)\otimes \mathcal L^{\otimes l}\right)
\\ &\to R^qf_*\left(\omega_X(D)\otimes \mathcal L^{\otimes k+l}\right)
\end{split}
\end{equation*} 
induced by $\otimes s$ is injective for every positive integer $l$. 
\end{itemize}
\end{thm}

Koll\'ar's original result in \cite{kollar1} is a special 
case of Theorem \ref{b-thm1.11}. 

\begin{say}[Koll\'ar's original statement]\label{b-say1.12}
If $X$ is a smooth projective variety with $D=0$ and 
$f\colon X\to Y$ is a projective surjective 
morphism onto a projective variety $Y$ in 
Theorem \ref{b-thm1.11} (i), then 
the strict support condition is nothing but 
Koll\'ar's torsion-freeness of $R^qf_*\omega_X$ 
(see \cite[Theorem 2.1 (i)]{kollar1}). 
We further assume that $Z$ is a point in Theorem \ref{b-thm1.11} (ii). 
Then we can recover Koll\'ar's vanishing theorem 
(see \cite[Theorem 2.1 (iii)]{kollar1}). 
If $X$ is a smooth projective variety, 
$D=0$, and $Y$ is a point, then Theorem \ref{b-thm1.11} (iii) 
coincides with Koll\'ar's original injectivity theorem 
(see \cite[Theorem 2.2]{kollar1}). 
Hence Theorem \ref{b-thm1.11} generalizes Koll\'ar's original 
statement in \cite{kollar1}. 
\end{say}

Our approach to Theorem \ref{b-thm1.11} 
in \cite{fujino-vanishing-ne}, which is 
completely different from 
the argument in \cite[Chapter 5]{fujino-foundations}, 
is very simple. 
By using a spectral sequence coming from 
Saito's theory of mixed Hodge modules (see Theorem 
\ref{b-thm2.3} below), 
we can reduce Theorem \ref{b-thm1.11} to 
a well-known simpler case due to Takegoshi (see Theorem 
\ref{b-thm2.1} below). Hence our proof of 
Theorem \ref{b-thm1.11} in \cite{fujino-vanishing-ne} 
uses the semisimplicity of polarizable Hodge modules. 
The advantage of this approach is to clarify the meaning 
of the strict support condition in Theorem 
\ref{b-thm1.11} (i). 
We note that the above semisimplicity comes 
from the semisimplicity 
of polarizable variations of pure Hodge structure 
since polarizable Hodge modules are uniquely 
determined by their generic variations of pure Hodge structure. 
We note that the reader can find an alternative 
approach to Theorem \ref{b-thm1.11}, 
which is free from Saito's theory of mixed Hodge modules and 
only depends on the semisimplicity of 
polarizable variations of pure Hodge structure, 
in \cite{fujino-fujisawa}.

\section{Sketch of Proof}\label{b-sec2}

In this section, we will briefly discuss how to prove the theorems 
in Section \ref{b-sec1}. As we have already 
explained in Section \ref{b-sec1}, 
we reduce the problem to a well-known simpler case 
due to Takegoshi by using a spectral sequence coming 
from the theory of mixed Hodge modules. 

The following theorem is a special case of Takegoshi's 
result (see \cite{takegoshi}). 
It is a complex analytic generalization of 
Koll\'ar's torsion-freeness and vanishing theorem. 

\begin{thm}[{see \cite{takegoshi}}]\label{b-thm2.1}
Let $f\colon X\to Y$ be a projective surjective morphism 
from a smooth irreducible complex analytic space $X$. 
Then $R^qf_*\omega_X$ is a torsion-free coherent sheaf on $Y$ 
for every $q$. 

Furthermore, let $\pi\colon Y\to Z$ be a projective morphism 
between complex analytic spaces and let $\mathcal A$ be a 
$\pi$-ample line bundle on $Y$. 
Then $R^p\pi_*\left(\mathcal A\otimes 
R^qf_*\omega_X\right)=0$ holds for 
every $p>0$ and every $q$. 
\end{thm}

Although Takegoshi's complex analytic approach 
to Koll\'ar's theorems in \cite{takegoshi} is interesting, 
we do not discuss it here since the statement of Theorem \ref{b-thm2.1} 
is sufficient for our purposes in this paper. 
For an alternative approach to Theorem \ref{b-thm2.1}, 
we recommend the reader to 
see \cite[Corollaries 1.2 and 1.5]{fujino-trans}. 

\begin{say}\label{b-say2.2}
Let $(X, D)$ be an analytic simple normal crossing pair such that 
$D$ is reduced. 
For any positive integer $k$, 
we put 
\begin{equation*}
X^{[k]}:=\left\{x\in X |\mult_x X\geq k\right\}^{\sim},  
\end{equation*}
where $Z^{\sim}$ denotes the normalization of $Z$. 
Then $X^{[k]}$ is the disjoint union of the intersections of $k$ irreducible 
components of $X$, and is smooth. 
We have a reduced simple normal crossing divisor 
$D^{[k]}\subset X^{[k]}$ defined by the pull-back of $D$ by the natural 
morphism 
$X^{[k]}\to X$. 
For any nonnegative integer $l$, 
we put 
\begin{equation*}
D^{[k, l]}:=\left\{\left. x\in X^{[k]} \right| \mult _x D^{[k]}\geq l\right\}^{\sim}. 
\end{equation*} 
We note that $D^{[k,0]}=X^{[k]}$ holds by definition. 
We also note that $\dim D^{[k, l]}=n+1-k-l$, where 
$n=\dim X$. 
In this situation, $W$ is a stratum of $(X, D)$ if and only if 
$W$ is the image of an irreducible component of 
$D^{[k, l]}$ for some $k>0$ and $l\geq 0$. 
\end{say}

One of the main ingredients of \cite{fujino-vanishing-ne} is 
the following result coming from Saito's theory of mixed 
Hodge modules (see \cite{saito1} and \cite{saito2}). 

\begin{thm}
[{\cite[Corollary 1 and 4.7.~Remark]{fujino-fujisawa-saito}}]\label{b-thm2.3}
Let $(X, D)$ be an analytic simple normal crossing pair 
with $\dim X=n$ such that 
$D$ is reduced and let $f\colon X\to Y$ 
be a projective morphism to a smooth 
complex analytic space $Y$.
Then there is the weight spectral sequence 
\begin{equation*}
\begin{split}
{}_FE^{-q, i+q}_1&=\bigoplus _{k+l=n+q+1}R^if_*\omega_{D^{[k, l]}/Y} 
\\&\Rightarrow R^if_*\omega_{X/Y}(D), 
\end{split}
\end{equation*} 
degenerating at $E_2$, and its $E_1$-differential 
$d_1$ splits so that the 
${}_FE^{-q, i+q}_2$ are direct factors of ${}_FE^{-q, i+q}_1$. 
\end{thm}

For the details of the weight spectral sequence of 
mixed Hodge modules on 
$Y$ necessary for Theorem \ref{b-thm2.3}, 
see \cite{fujino-fujisawa-saito}. 
Once we know Theorems \ref{b-thm2.1} and \ref{b-thm2.3}, 
it is easy to prove Theorem \ref{b-thm1.11}. 
Here, we only prove Theorem \ref{b-thm1.11} (i) and (ii). 

\begin{proof}[Proof of Theorem \ref{b-thm1.11}]
First we prove (i). 
Since the problem is local, we may assume that 
$Y$ is a closed analytic subspace of a polydisc $\Delta^m$. 
By replacing $Y$ with $\Delta^m$, 
we may further assume that $Y$ itself is a polydisc. 
In this case, we can use Theorem \ref{b-thm2.3}. 
We note that $\omega_Y\simeq \mathcal O_Y$ holds. 
By Theorem \ref{b-thm2.1}, 
\begin{equation*}
{}_FE^{-q, i+q}_1\simeq \bigoplus _{k+l=n+q+1}R^if_*\omega_{D^{[k, l]}}
\end{equation*}
satisfies the strict support condition, that is, 
every associated subvariety of 
\begin{equation*}
{}_FE^{-q, i+q}_1\simeq \bigoplus _{k+l=n+q+1}R^if_*\omega_{D^{[k, l]}}
\end{equation*} 
is the $f$-image of some stratum of $(X, D)$.   
By Theorem \ref{b-thm2.3}, the associated subvariety 
of ${}_FE^{-q, i+q}_2={}_FE^{-q, i+q}_\infty$ 
is the $f$-image of some stratum of $(X, D)$. 
This implies that $R^qf_*\omega_X(D)$ satisfies the 
desired strict support condition. 
Next, we treat (ii). 
We may assume that $Z$ is a polydisc and $Y$ is a closed analytic 
subspace of $Z\times \mathbb P^n$. 
By applying Theorem \ref{b-thm2.3} to 
$f\colon X\to Y\hookrightarrow Z\times \mathbb P^n$, 
we obtain the following spectral sequence 
\begin{equation*}
E^{-q, i+q}_1=\bigoplus _{k+l=n+q+1}R^if_*\omega_{D^{[k, l]}} 
\Rightarrow R^if_*\omega_X(D)
\end{equation*}
which degenerates at $E_2$ such that its $E_1$-differential $d_1$ splits. 
By Theorem \ref{b-thm2.1}, 
we obtain 
\begin{equation*}
R^p\pi_*\left(\mathcal A\otimes E^{-q, i+q}_1\right)=0
\end{equation*} 
for every 
$p>0$. Since the $E^{-q, i+q}_2=E^{-q, i+q}_\infty$ are direct factors of 
$E^{-q, i+q}_1$, we have 
\begin{equation*}
R^p\pi_*\left(\mathcal A\otimes E^{-q, i+q}_2\right)=0
\end{equation*}  
for every $p>0$. 
This implies that 
\begin{equation*} 
R^p\pi_*\left(\mathcal A\otimes R^qf_*\omega_X(D)\right)=0
\end{equation*}  
holds for every $p>0$. 
This is what we wanted. 
\end{proof}

Since the proof of Theorem \ref{b-thm1.6} in 
\cite{fujino-vanishing-ne} is somewhat technical, 
we only give a sketch here. 
For the details, see \cite[Sections 4 and 5]{fujino-vanishing-ne}. 

\begin{proof}[Sketch of Proof of Theorem \ref{b-thm1.6}]
For (i), we take an arbitrary point $y\in Y$. 
We may shrink $Y$ around $y$ suitably. 
By perturbing $\Delta$, we may further assume that 
$\Delta$ is a $\mathbb Q$-divisor. 
Then we reduce the problem to the case where 
$(X, \Delta)$ is an analytic globally embedded simple 
normal crossing pair (see \cite[Lemma 5.1]{fujino-vanishing-ne}). 
By repeatedly taking suitable 
finite covers, we can reduce the problem to 
Theorem \ref{b-thm1.11} (i). 
For (ii), we can shrink $Z$ suitably and may assume 
that $\Delta$ is a $\mathbb Q$-divisor. 
As for (i), we make $(X, \Delta)$ an 
analytic globally embedded simple 
normal crossing pair (see \cite[Lemma 5.1]{fujino-vanishing-ne}) 
and repeatedly take suitable finite 
covers. 
Then we see that (ii) follows from Theorem \ref{b-thm1.11} (ii). 
\end{proof}

Theorem \ref{b-thm1.9} follows from Theorem \ref{b-thm1.6}. 

\begin{proof}[Sketch of Proof of Theorem \ref{b-thm1.9}]
We need no new ideas for the proof of Theorem \ref{b-thm1.9}. 
The proof of \cite[Theorem 5.7.3]{fujino-foundations} 
can work with some suitable modifications. 
Theorem \ref{b-thm1.9} can be seen as a corollary of 
Theorem \ref{b-thm1.6}. 
\end{proof}
 
\begin{say}[Traditional approach versus new approach]\label{b-say2.4}
Here we briefly review the main difference 
between the approach adopted in 
\cite[Chapter 5]{fujino-foundations} and the 
one explained in this paper. Note that in 
\cite[Chapter 5]{fujino-foundations} 
everything is assumed to be {\em{algebraic}}. 
In \cite[Chapter 5]{fujino-foundations}, 
we first establish an injectivity theorem 
by using the theory of mixed Hodge 
structures on cohomology with compact 
support (see \cite[Sections 5.4 and 5.5]{fujino-foundations}). 
By this Hodge theoretic injectivity theorem, we 
can check that the injectivity theorem 
like Theorem \ref{b-thm1.11} (iii) holds true in the 
algebraic setting (see, for example, 
\cite[Theorem 5.6.1]{fujino-foundations}). 
It is well known that we can quickly recover 
Theorem \ref{b-thm1.11} (i) and (ii) once we 
obtain Theorem \ref{b-thm1.11} (iii). 
When $X$ is algebraic, we can always take a 
compactification $\overline X$ of $X$ and 
use the mixed Hodge structures on $\overline X$. 
On the other hand, if $X$ is not algebraic, 
then we can not always compactify $X$. Thus, 
the argument in \cite[Chapter 5]{fujino-foundations} 
does not work when $X$ is not algebraic. 
Anyway, a key result in the approach 
of \cite[Chapter 5]{fujino-foundations} is the 
injectivity theorem 
coming from the theory of 
mixed Hodge structures on cohomology with compact support. 
In the approach explained in this paper, injectivity 
theorems do not play an important role. 
\end{say}

\begin{ack}\label{b-ack}
The author thanks Professor Taro Fujisawa very much. 
He was partially 
supported by JSPS KAKENHI Grant Numbers 
JP19H01787, JP20H00111, JP21H00974, JP21H04994. 
\end{ack}

%%%%%%%%%%%%%%%%%

\end{document}